\documentclass[11pt]{amsart}
\usepackage{graphicx,amssymb}
\usepackage{a4wide}

\theoremstyle{plain}
\newtheorem{thm}{Theorem}[section]
\newtheorem{lem}[thm]{Lemma}
\newtheorem{prop}[thm]{Proposition}
\newtheorem{cor}[thm]{Corollary}

\theoremstyle{definition}
\newtheorem{defn}[thm]{Definition}
\newtheorem{exmp}[thm]{Example}

\newenvironment{pf}{\begin{proof}}{\end{proof}}

\def\Z{\mathbb Z}
\def\Q{\mathbb Q}
\def\N{\mathbb N}
\def\c{\mathbf c}
\def\d{\mathbf d}
\def\D{\mathcal D}
\def\DD{\mathbf D}
\def\wedgeL{\wedge_L}
\def\wedgeR{\wedge_R}
\def\veeL{\vee_{L}}
\def\veeR{\vee_{R}}
\def\Aut{\operatorname{Aut}}
\let\le\leqslant
\let\ge\geqslant
\def\len{\operatorname{len}}
\def\infs{\inf{\!}_ s}
\def\sups{\sup{\!}_s}
\def\lens{\len{\!}_s}

\begin{document}

\title{Garside groups are strongly translation discrete}

\author{Sang Jin Lee}
\thanks{This paper was supported by Konkuk University in 2005}
\address{Department of Mathematics, Konkuk University, Gwangjin-gu, Seoul 143-701, Korea}
\email{sangjin@konkuk.ac.kr}
\date{\today}

\begin{abstract}
In this paper we show that all Garside groups are strongly
translation discrete, that is, the translation numbers of
non-torsion elements are strictly positive and for any real number
$r$ there are only finitely many conjugacy classes of elements whose
translation numbers are less than or equal to $r$. It is a
consequence of the inequality
``$\infs(g)\le\frac{\infs(g^n)}n<\infs(g)+1$'' for a positive
integer $n$ and an element $g$ of a Garside group $G$, where $\infs$
denotes the maximal infimum for the conjugacy class. We prove the
inequality by studying the semidirect product $G(n)=\Z\ltimes G^n$
of the infinite cyclic group $\Z$ and the cartesian product $G^n$ of
a Garside group $G$, which turns out to be a Garside group. We also
show that the root problem in a Garside group $G$ can be reduced to
a conjugacy problem in $G(n)$, hence the root problem is solvable
for Garside groups.

\end{abstract}

\maketitle

\section{Introduction}

The notion of translation numbers of finitely generated groups
was first introduced by Gersten and Short in~\cite{GS91}.
This concept comes from the action of the fundamental group of a compact
Riemannian manifold of non-positive curvature
on the universal cover of this manifold.

\begin{defn}
Let $G$ be a finitely generated group and $X$
a finite set of semigroup generators for $G$.
The \emph{translation number} with respect to $X$ of
a non-torsion element $g\in G$ is defined by
$$t_X(g)=\lim_{n\to \infty}\frac{|g^n|{}_X}n,$$
where $|\cdot|_X$ denotes the shortest word length
in the alphabet $X$.
\end{defn}

For $g\in G$, let $[g]$ denote the conjugacy class of $g$ in $G$.
As noted in~\cite{GS91}, the translation number is
constant on each conjugacy class in $G$.
The following definitions were suggested
by Conner~\cite{Con00} and Kapovich~\cite{Kap97}.

\begin{defn}
A finitely generated group $G$ is said to be
\begin{enumerate}
\item \emph{translation separable} (or \emph{translation proper})
if for some (and hence for any) finite set $X$
of semigroup generators for $G$
the translation numbers of non-torsion elements
are strictly positive;
\item \emph{translation discrete} if it is translation separable
and for some (and hence for any) finite set $X$
of semigroup generators for $G$
the set $t_X(G)$ has 0 as an isolated point;
\item \emph{strongly translation discrete} if it is translation separable
and for some (and hence for any) finite set $X$
of semigroup generators for $G$
and for any real number $r$
the number of conjugacy classes $[g]$ in $G$
with $t_X(g)\le r$ is finite.
\end{enumerate}
\end{defn}

There are good properties of groups which have some kind of
translation discreteness.
Gersten and Short~\cite{GS91} proved
that every finitely generated nilpotent subgroup
of a translation separable group is virtually abelian.
Conner~\cite{Con00} proved
that a translation separable solvable group
of finite virtual cohomological dimension is metabelian by finite
and that every solvable subgroup of finite virtual cohomological dimension
in a translation discrete group is a finite extension of $\Z^m$.

Therefore, it would be interesting to study translation discreteness
of groups. Biautomatic groups are translation separable by Gersten
and Short~\cite{GS91}. Word hyperbolic groups are strongly
translation discrete. Moreover, the translation numbers in a word
hyperbolic group are rational and have bounded denominators.
This follows from a claim of Gromov~\cite[5.2C]{Gr87} and was accurately
proved by Swenson~\cite[Corollary of Theorem 13]{Swe95}.
Kapovich~\cite{Kap97} proved that C(4)-T(4)-P, C(3)-T(6)-P and
C(6)-P small cancellation groups are strongly translation discrete.
(These groups are biautomatic but not necessarily word hyperbolic by
Gersten and Short~\cite{GS90}.) Bestvina~\cite{Be99} showed that the
Artin groups of finite type are translation discrete.

The class of Garside groups provides a lattice-theoretic
generalization of the braid groups and the Artin groups of finite type.
Garside groups have been studied actively for the last ten years.
Recently, Charney, Meier and Whittlesey~\cite{CMW04}
proved that a Garside group is translation discrete
under the condition that it has a tame Garside element,
generalizing the work of Bestvina.
See \S\ref{sec:trans} for the definition of tameness.
In this paper, we show

\par\medskip
\noindent\textbf{Main Theorem}(Theorem~\ref{thm:main})\hskip 1em\em
All Garside groups are strongly translation discrete.
\normalfont\medskip

We remark that we don't need the existence of a tame Garside element.
Moreover, even when restricted to the Artin groups of finite type,
our result is stronger than that of Bestvina in~\cite{Be99}.

Our approach is different from that of Bestvina
(and so that of Charney, Meier and Whittlesey).
We show that if $G$ is a Garside group,
then so is the semidirect product $G(n)=\Z\ltimes G^n$
of the infinite cyclic group $\Z$ and the cartesian product $G^n$,
where the action of $\Z$ on $G^n$ is a coordinate change.
(It is a special case of the crossed product of Picantin~\cite{Pic01a}.)
Studying the relationship between the Garside structure of $G(n)$
and that of $G$, we show that
$$\infs(g)\le\frac{\infs(g^n)}n<\infs(g)+1$$
for all $n\ge 1$, where $\infs$ denotes the maximal infimum
for the conjugacy class.
The main theorem follows from it.

\smallskip
Another result is a new solution to the root problem
in Garside groups.
The root problem in a group $G$ is
to solve the equation $x^n=g$ in $G$ for $x$,
given $g\in G$ and $n\ge 2$.
Sty\v shnev~\cite{Sty78} showed that the problem is solvable
for the braid groups, and Sibert~\cite{Sib02}
obtained the same result for the Garside groups
under the assumption that the positive conjugacy classes
are all finite.
See \S\ref{sec:root} for the exact description of his assumption.
Roughly, Sty\v shnev and Sibert showed the following.

\begin{quote}
Let $G$ be a Garside group and $G^+$ its positive monoid.
Given $g\in G$ and $n\ge 2$, there exists a finite set
$C\subset G^+$
such that $x^n=g$ has a solution in $G$ if and only if $x^n=c$ has
a solution in $G^+$ for some $c\in C$.
\end{quote}

Note that we can solve the equation $x^n=c$ in $G^+$
using a sort of exhaustive search.
In this paper we show that
the root problem in a Garside group $G$
can be reduced to a conjugacy problem in $G(n)$.
Our solution does not need any extra condition.

\section{Garside monoids and groups}

In the late sixties, Garside~\cite{Gar69} solved the word and conjugacy
problems in the braid groups studying the positive braid monoids.
His theory has been generalized and
improved by several
mathematicians~\cite{BS72,Thu92,EM94,BKL98,DP99,BKL01,FG03,Pic01b,Geb03}.
Recently, Dehornoy and Paris~\cite{DP99} introduced
the notion of Garside groups,
which is a lattice-theoretic generalization
of the braid groups and the Artin groups of finite type.
Here we collect relevant information about
the word and conjugacy problems in Garside groups.
See~\cite{DP99,Deh02,FG03,Pic01b,Geb03} for details.

Let $M$ be a monoid.
Let \emph{atoms} (or `indivisible elements') be the elements $a\in M$
such that $a\ne 1$ and if $a=bc$ then either $b=1$ or $c=1$.
For $a\in M$, let $\Vert a\Vert$ be the supremum
of the lengths of all expressions of
$a$ in terms of atoms. The monoid $M$ is said to be \emph{atomic}
if it is generated by its atoms and $\Vert a\Vert<\infty$ for any $a\in M$.
In an atomic monoid $M$, there are partial orders $\le_L$ and $\le_R$:
$a\le_L b$ if $ac=b$ for some $c\in M$;
$a\le_R b$ if $ca=b$ for some $c\in M$.
For $a\in M$, let $L(a)=\{b\in M:b\le_L a\}$ and $R(a)=\{b\in M: b\le_R a\}$.

\begin{defn}
An atomic monoid $M$ is called a \emph{Garside monoid} if
\begin{enumerate}
\item $M$ is finitely generated;
\item $M$ is left and right cancellative;
\item $(M,\le_L)$ and $(M,\le_R)$ are lattices;
\item there exists an element $\Delta$, called a
\emph{Garside element}, such that
$L(\Delta)$ and $R(\Delta)$ are the same and they form
a set of generators for $M$.
\end{enumerate}
\end{defn}

The elements of $L(\Delta)$ are called the \emph{simple elements}.
Let $\D$ be the set of simple elements, that is, $\D=L(\Delta)$.
Let $\wedgeL$ and $\veeL$ (respectively, $\wedgeR$ and $\veeR$)
denote the gcd and lcm with respect to $\le_L$ (respectively, $\le_R$).

Garside monoids satisfy Ore's conditions,
and thus embed in their groups of fractions.
A \emph{Garside group} is defined to be the group of fractions
of a Garside monoid.
When $M$ is a Garside monoid and $G$ the group of fractions of $M$,
we identify the elements of $M$ and their images in $G$
and call them the \emph{positive elements} of $G$.
$M$ is called the \emph{positive monoid} of $G$,
often denoted $G^+$.
The partial orders $\le_L$ and $\le_R$, and thus the lattice structures
in the positive monoid $G^+$ can be extended
to the Garside group $G$ as follows:
$g\le_L h$ (respectively, $g\le_R h$) for $g,h\in G$
if $gc=h$ (respectively, $cg=h$) for some $c\in G^+$.
Let $\tau\colon G\to G$ be the inner automorphism of $G$ defined by
$\tau(g)=\Delta^{-1} g\Delta$.

\begin{thm}
Let $G$ be a Garside group, $G^+$ its positive monoid
and $\Delta$ a Garside element.
\begin{enumerate}
\item $\tau(G^+)=G^+$.
\item There is an integer $e$ such that $\Delta^e$ is central in $G$.
\item For $g\in G$, there are integers $r$ and $s$ such that
$\Delta^r\le_L g\le_L\Delta^s$.
\item For $g\in G$, there are unique simple elements
$s_1,\ldots,s_k\in \D\setminus\{1,\Delta\}$ such that
$$g=\Delta^r s_1\cdots s_k,$$
and $(s_is_{i+1}\cdots s_k)\wedgeL \Delta=s_i$ for $i=1,\ldots,k$.
\end{enumerate}
\end{thm}

By (3) the following invariants are well-defined:
$\inf(g)=\max\{r\in\Z:\Delta^r\le_L g\}$;
$\sup(g)=\min\{s\in\Z:g\le_L \Delta^s\}$;
$\len(g)=\sup(g)-\inf(g)$.
The expression in (4) is called the \emph{normal form} of $g$.
If $\Delta^r s_1\cdots s_k$ is the normal form of $g$,
then $\inf(g)=r$ and\/ $\sup(g)=r+k$.

We remark that $\tau$ and the invariants $\inf(g)$, $\sup(g)$
and $\len(g)$ depend on the choice of the Garside element.
(In \S3, we discuss non-uniqueness of Garside elements.)
Throughout the paper, when we mention a Garside monoid/group
without specifying its Garside element,
we always assume that a Garside element is chosen and fixed.

\begin{defn}
For $g\in G$, $[g]$ denotes its conjugacy class.
We define $\infs(g)=\max\{\inf(h):h\in [g]\}$,
$\sups(g)=\min\{\sup(h):h\in [g]\}$ and $\lens(g)=\sups(g)-\infs(g)$.
\end{defn}

\begin{defn}
For $g\in G$, the set $[g]^S=\{h\in [g]:
\inf(h)=\infs(g),\ \sup(h)=\sups(g)\}$
is called the \emph{super summit set} of $g$.
\end{defn}

\begin{defn}
Let $\Delta^r s_1\cdots s_k$ be the normal form of\/ $g\in G$.
Then $\c(g)=\Delta^r s_2\cdots s_k \tau^{-r}(s_1)$
and $\d(g)=\Delta^r\tau^r(s_k)s_1\cdots s_{k-1}$
 are called the \emph{cycling}
and \emph{decycling} of $g$.
\end{defn}

\begin{thm}\label{thm:SSS}
Let $G$ be a Garside group and $g\in G$.
\begin{enumerate}
\item If\/ $\c^k(g)=g$ for some $k\ge 1$, then $\inf(g)=\infs(g)$.
\item If\/ $\d^k(g)=g$ for some $k\ge 1$, then $\sup(g)=\sups(g)$.
\item If\/ $h\in [g]^S$, then $\c(h),\d(h),\tau(h)\in [g]^S$.
\item If\/ $h\in [g]^S$, then
    $\tau(\c(h))=\c(\tau(h))$ and $\tau(\d(h))=\d(\tau(h))$.
\item $\c^k(\d^l(g))\in[g]^S$ for some $k,l\ge 0$.
\item $[g]^S$ is finite and non-empty.
\item For $h,h'\in[g]^S$, there is a finite sequence
$h=h_0\to h_1\to\cdots\to h_m=h'$
such that for $i=1,\ldots,m$, $h_i\in[g]^S$ and $h_{i}=s_i^{-1}h_{i-1}s_i$
for some $s_i\in\D$.
\end{enumerate}
\end{thm}

Recently, Franco and Gonz\'alez-Meneses~\cite{FG03} improved the algorithm
to generate the super summit set $[g]^S$ from an element of $[g]^S$.
Gebhardt~\cite{Geb03} defined the \emph{ultra summit set}
of $g\in G$ as
$[g]^U=\{h\in [g]^S:\mbox{$\c^k(h)=h$ for some $k>0$}\}$
and showed that it satisfies most important properties
of the super summit set. (Theorem~\ref{thm:SSS} holds when we replace
the super summit set $[g]^S$ with the ultra summit set $[g]^U$.)

We will often use the following lemma later.
It is easy to prove and so we omit the proof.
\begin{lem}\label{thm:dual}
Let $G$ be a Garside group and $g\in G$.
$$\d(g)=(\tau^{-1}(\c(g^{-1})))^{-1},\quad
\sup(g)=-\inf(g^{-1}),\quad
\sups(g)=-\infs(g^{-1}).$$
\end{lem}

\section{Minimal Garside element}

In a Garside monoid, Garside elements are not unique.
For example, if $\Delta$ is a Garside element,
then $\Delta^m$ is also a Garside element for any $m\ge 1$.
It is known by Dehornoy that there exists a unique Garside element
whose $\Vert\cdot\Vert$-norm is minimal.
See~\cite[the discussion after Definition 1.11 in page 275]{Deh02}.
In this section, we give an another proof.

\begin{lem}\label{thm:lemma1}
Let $G^+$ be a Garside monoid and $a,b\in G^+$.
If\/ $L(a)=R(b)$, then $a=b$.
\end{lem}

\begin{pf}
Since $a\in L(a)=R(b)$ and $b\in R(b)=L(a)$, $a\le_R b$ and $b\le_L a$.
Therefore, $b=c_1a$ and $a=bc_2$ for some $c_1,c_2\in G^+$,
and thus  $b=c_1a=c_1bc_2$. Since
$$
\Vert b\Vert
=\Vert c_1bc_2\Vert\ge \Vert c_1\Vert + \Vert b\Vert + \Vert c_2\Vert
\ge \Vert b\Vert,
$$
$\Vert c_1\Vert=\Vert c_2\Vert=0$.
Hence, $c_1=c_2=1$, and thus $a=b$.
\qed\end{pf}

\begin{lem}\label{thm:GarsideWedge}
If\/ $\Delta_1$ and $\Delta_2$ are Garside elements of a Garside monoid\/ $G^+$,
then $\Delta_1\wedgeL\Delta_2=\Delta_1\wedgeR\Delta_2$ and
it is a  Garside element of\/ $G^+$.
\end{lem}

\begin{pf}
Note that for $a,b\in G^+$,
$L(a\wedgeL b)=L(a)\cap L(b)$ and $R(a\wedgeR b)=R(a)\cap R(b)$.
Since $L(\Delta_i)=R(\Delta_i)$ for $i=1,2$,
$$
L(\Delta_1\wedgeL \Delta_2)
=L(\Delta_1)\cap L(\Delta_2)
=R(\Delta_1)\cap R(\Delta_2)
=R(\Delta_1\wedgeR\Delta_2).
$$
By Lemma~\ref{thm:lemma1},
$\Delta_1\wedgeL\Delta_2=\Delta_1\wedgeR\Delta_2$.
Since the generating subsets of an atomic monoid are exactly
those subsets that contain all atoms,
the intersection of any collection of generating subsets
is a generating subset.
Therefore $L(\Delta_1\wedgeL\Delta_2)=L(\Delta_1)\cap L(\Delta_2)$
is a generating subset of $G^+$, and thus
$\Delta_1\wedgeL\Delta_2$ is a Garside element.
\qed\end{pf}

\begin{lem}\label{thm:minimalGarside}
In a Garside monoid\/ $G^+$, there is a unique Garside element
that is minimal with respect to both $\le_L$ and $\le_R$
among all Garside elements of\/ $G^+$.
\end{lem}

\begin{pf}
Since a Garside monoid is finitely generated,
$G^+$ has countably many elements,
and thus countably many Garside elements,
say $\{\Delta_j\}_{j\ge 1}$.
Let $\Delta^{(i)}=\Delta_1\wedge_L\Delta_2\wedge_L\cdots\wedge_L\Delta_i$.
Since $\Delta^{(i+1)}\le_L\Delta^{(i)}$ for all $i$,
the atomicity of $G^+$ implies $\Delta^{(m)}=\Delta^{(m+1)}=\cdots$
for some $m\ge 1$.
By Lemma~\ref{thm:GarsideWedge}, $\Delta^{(m)}$ is a Garside element.
Since it is both the left and right gcd of all Garside elements,
it is minimal with respect to both $\le_L$ and $\le_R$
among all Garside elements of $G^+$.
Uniqueness is obvious.
\qed\end{pf}

\begin{defn}
The unique Garside element of Lemma~\ref{thm:minimalGarside} is
called the \emph{minimal Garside element} of the Garside monoid $G^+$.
\end{defn}

\begin{lem}\label{thm:AutoProperty}
Let $\phi$ be a monoid automorphism of a Garside monoid\/ $G^+$.
\begin{enumerate}
\item $\phi$ permutes the atoms of\/ $G^+$.
\item For $a,b\in G^+$,
$a\le_L b$ if and only if\/ $\phi(a)\le_L\phi(b)$.
The same is true for $\le_R$.
\item For $a,b\in G^+$, $\phi(a\wedgeL b)=\phi(a)\wedgeL\phi(b)$.
The same is true for $\veeL$, $\wedgeR$ and $\veeR$.
\item If\/ $\phi(\Delta)=\Delta$ for a Garside element $\Delta$,
then $\phi$ permutes the simple elements.
\item If\/ $\Delta$ is the minimal Garside element of\/ $G^+$,
then $\phi(\Delta)=\Delta$.
\end{enumerate}
\end{lem}

\begin{pf}
(1), (2) and (3) are easy. See~\cite[Lemma 9.1, p.\ 599]{DP99}.
(4) is immediate from (2). Let us prove (5).
Because the automorphism $\phi$ respects the partial orders
$\le_L$ and $\le_R$, $\phi(\Delta)$ is the minimal Garside element of $G^+$,
hence $\phi(\Delta)$ is equal to $\Delta$ by the uniqueness of
the minimal Garside element.
\qed\end{pf}

The following theorem summarizes the discussions in this section.

\begin{thm}
In a Garside monoid\/ $G^+$, there exists a unique Garside element $\Delta$
that is minimal with respect to both $\le_L$ and $\le_R$
among all Garside elements of\/ $G^+$.
Furthermore, for any monoid automorphism $\phi$ of\/ $G^+$,
$\phi(\Delta)=\Delta$.
\end{thm}

\section{Semidirect product of Garside monoids}

Let $G$ and $H$ be Garside groups.
Let $G$ (and hence $G^+$) act on $H^+$ on the right
via a homomorphism $\rho\colon G\to\Aut(H^+)$,
where $\Aut(H^+)$ is the group of automorphisms of the monoid $H^+$.
We use the notation $b^g$ to denote the action of $g\in G$ on
$b\in H^+$ via $\rho$, that is, $b^g=\rho(g)(b)$.
Recall that the semidirect product $G^+\ltimes_\rho H^+$ has
the underlying set $\{(a,b):a\in G^+, b\in H^+\}$ and the product
$$
(a_1,b_1)(a_2,b_2)=(a_1a_2,b_1^{a_2}b_2),\qquad
\mbox{$a_1,a_2\in G^+$ and $b_1,b_2\in H^+$.}
$$
Let $G\ltimes_\rho H$ denote the semidirect product of
the Garside groups $G$ and $H$ via the homomorphism
$G\to\Aut(H)$ that is the composition of $\rho$
with the canonical monomorphism $\Aut(H^+)\to \Aut(H)$.
If there is no confusion, we omit $\rho$
from the notation of semidirect products.
The main goal of this section is to prove the following theorem.

\begin{thm}\label{thm:Semidirect}
If\/ $G$ and $H$ are Garside groups and $\rho:G\to\Aut(H^+)$ is a homomorphism,
then $G^+\ltimes_\rho H^+$ is a Garside monoid and its group of
fractions is $G\ltimes_\rho H$.
Moreover, if\/ $\Delta_G$ and $\Delta_H$ are Garside elements of\/ $G$ and $H$,
and\/ $\Delta_H$ is fixed under the action of\/ $G$,
then $(\Delta_G,\Delta_H)$ is a Garside element of\/ $G^+\ltimes_\rho H^+$
and the set of simple elements is $\{(a,b):a\le_L \Delta_G, b\le_L \Delta_H\}$.
\end{thm}

In fact, Theorem~\ref{thm:Semidirect} is already known by Picantin.
He showed that the crossed product of Garside monoids
is a Garside monoid~\cite[Proposition 3.12]{Pic01a}.
Our semidirect product is a special case of his crossed product.
Because his construction is more general than ours,
the proof in his paper is longer and more complicated.
So we include the proof.
Theorem~\ref{thm:Semidirect} is a consequence of the following proposition.

\begin{prop}\label{thm:SemidirectProperties}
Let $G^+\ltimes H^+$ and $G\ltimes H$ be the semidirect products of
Garside monoids and groups via a homomorphism
$\rho:G\to\Aut(H^+)$.
Let $(a_1,b_1),(a_2,b_2)\in G^+\ltimes H^+$.
\begin{enumerate}
\item $G^+\ltimes H^+$ is (left and right) cancellative.

\item $G^+\ltimes H^+$ is atomic.

\item $(a_1,b_1)\le_L(a_2,b_2)$ if and only if $a_1\le_L a_2$ and
$b_1^{a_1^{-1}}\le_L b_2^{a_2^{-1}}$.

\item $(a_1,b_1)\le_R(a_2,b_2)$ if and only if $a_1\le_R a_2$ and
$b_1\le_R b_2$.

\item $(a_1,b_1)\wedgeL(a_2,b_2)=(a_1\wedgeL a_2,
(b_1^{{a_1}^{-1}}\wedgeL b_2^{{a_2}^{-1}})^{a_1\wedgeL a_2})$.

\item $(a_1,b_1)\veeL (a_2,b_2)=(a_1\veeL a_2,
(b_1^{{a_1}^{-1}}\veeL b_2^{{a_2}^{-1}})^{a_1\veeL a_2})$.

\item $(a_1,b_1)\wedgeR(a_2,b_2)=(a_1\wedgeR a_2, b_1\wedgeR
b_2)$.

\item $(a_1,b_1)\veeR(a_2,b_2)=(a_1\veeR a_2, b_1\veeR b_2)$.

\item Let $\Delta_G$ and $\Delta_H$ be Garside elements of\/ $G$
and $H$. If\/ $\Delta_H$ is fixed under the action of\/ $G$, then
$(\Delta_G,\Delta_H)$ is a Garside element of\/ $G^+\ltimes H^+$
and the set of simple elements is $\{(a,b):a\le_L \Delta_G, b\le_L
\Delta_H\}$.

\item $G^+\ltimes H^+$ is a Garside monoid and its group of
fractions is $G\ltimes H$.
\end{enumerate}
\end{prop}

\begin{pf}
(1) is obvious.

\smallskip
(2)\ \
The atoms of $G^+\ltimes H^+$ are of the form either
$(a,1)$ for an atom $a\in G^+$ or $(1,b)$ for an atom $b\in H^+$.
Therefore, $G^+\ltimes H^+$ is generated by its atoms.
Since $\|(a,b)\|=\|a\|+\|b\|<\infty$ for $(a,b)\in G^+\ltimes H^+$,
$G^+\ltimes H^+$ is atomic.

\smallskip
(3)\ \
If $(a_1,b_1)\le_L (a_2,b_2)$,
then $(a_2,b_2)=(a_1,b_1)(a_3,b_3)=(a_1a_3,b_1^{a_3}b_3)$
for some $(a_3,b_3)\in G^+\ltimes H^+$.
Since $a_2=a_1a_3$, $a_1\le_L a_2$.
Since $b_2=b_1^{a_3}b_3=b_1^{{a_1}^{-1}a_2}b_3$,
$b_1^{{a_1}^{-1}a_2}\le_L b_2$,
and thus $b_1^{{a_1}^{-1}}\le_L b_2^{{a_2}^{-1}}$.
Conversely, if $a_1\le_L a_2$ and $b_1^{{a_1}^{-1}}\le_L b_2^{{a_2}^{-1}}$,
then $a_2=a_1a_3$ and $b_2^{{a_2}^{-1}}=b_1^{{a_1}^{-1}}b_3$
for some $(a_3,b_3)\in G^+\ltimes H^+$.
Since $(a_1,b_1)(a_3,b_3^{a_2})=(a_1a_3,b_1^{a_3}b_3^{a_2})
=(a_2,b_1^{{a_1}^{-1}a_2}b_3^{a_2})
=(a_2,(b_1^{{a_1}^{-1}}b_3)^{a_2})
=(a_2,b_2)$,
$(a_1,b_1)\le_L (a_2,b_2)$.

\smallskip
(4) can be proved similarly to (3).

\smallskip
(5)\  $\begin{array}[t]{rcl}
\lefteqn{(a_3,b_3)\le_L (a_1,b_1)\mbox{ and }(a_3,b_3)\le_L (a_2,b_2)}\\
&&\Leftrightarrow a_3\le_L a_1,\quad a_3\le_L a_2,\quad
b_3^{{a_3}^{-1}}\le_L b_1^{{a_1}^{-1}},\quad b_3^{{a_3}^{-1}}\le_L b_2^{{a_2}^{-1}}
\mbox{}\qquad\mbox{by (3)}\\
&&\Leftrightarrow a_3\le_L (a_1\wedgeL a_2)\mbox{ and }
b_3^{{a_3}^{-1}}\le_L (b_1^{{a_1}^{-1}}\wedgeL b_2^{{a_2}^{-1}})
\mbox{}\qquad\mbox{by the definition of $\wedgeL$}\\
&&\Leftrightarrow (a_3,b_3)\le_L (a_1\wedgeL a_2,
(b_1^{{a_1}^{-1}}\wedgeL b_2^{{a_2}^{-1}})^{(a_1\wedgeL a_2)})
\mbox{}\qquad\mbox{by (3).}
\end{array}$

\smallskip
(6), (7) and (8) can be proved similarly to (5).

\smallskip
(9)\ \
$L((\Delta_G,\Delta_H))=R((\Delta_G,\Delta_H))$ because
$$\begin{array}[t]{l}
(a,b)\le_L (\Delta_G,\Delta_H)\\
\mbox{}\quad \Leftrightarrow a\le_L \Delta_G,\ \
         b^{a^{-1}}\le_L(\Delta_H)^{{\Delta_G}^{-1}}
\mbox{}\qquad\mbox{by (3)}\\
\mbox{}\quad \Leftrightarrow a\le_L \Delta_G,\ \ b\le_L \Delta_H
\mbox{}\qquad\mbox{since $\Delta_H$ is fixed under the action of $G$}\\
\mbox{}\quad \Leftrightarrow a\le_R \Delta_G,\ \ b\le_R\Delta_H
\mbox{}\qquad\mbox{by definition of Garside elements}\\
\mbox{}\quad \Leftrightarrow (a,b)\le_R(\Delta_G,\Delta_H)
\mbox{}\qquad\mbox{by (4)}.
\end{array}
$$
In particular, $R((\Delta_G,\Delta_H))
=\{(a,b): a\le_R\Delta_G, b\le_R\Delta_H\}$
contains all the atoms of $G^+\ltimes H^+$,
and thus it generates $G^+\ltimes H^+$.
Hence $(\Delta_G,\Delta_H)$ is a Garside element and the set of simple elements is
as in the statement.

\smallskip
(10)\ \
$G^+\ltimes H^+$ is finitely generated
since $G^+$ and $H^+$ are finitely generated;
cancellative by (1);
atomic by (2).
$(G^+\ltimes H^+,\le_L)$ and $(G^+\ltimes H^+,\le_R)$ are lattices
by (5), (6), (7) and (8).
Since, at least, the minimal Garside element of $H$ is fixed under the action
of $G$, $G^+\ltimes H^+$ has a Garside element by (9).
Hence $G^+\ltimes H^+$ is a Garside monoid.
It is easy to see that
$G\ltimes H$ is the group of fractions of $G^+\ltimes H^+$.
\qed\end{pf}

We remark that our semidirect product is different from that
of Crisp and Paris in~\cite{CP05}.
They studied the semidirect product of the form
$(H*\cdots *H)\rtimes_\rho B_n$, where $B_n$ is the $n$-braid group
and $\rho:B_n\to\Aut(H\ast\cdots\ast H)$ is an Artin type representation.

\begin{exmp}
(1) For a Garside group $G$, the cartesian product $G^n=G\times\cdots \times G$
is a Garside group with a Garside element
$(\Delta,\ldots,\Delta)$, where $\Delta$ is a Garside element of $G$.

(2) The $n$-braid group $B_n$ acts on $G^n$ by
$(g_1,\ldots,g_n)^\alpha=(g_{\theta^{-1}(1)},\ldots,g_{\theta^{-1}(n)})$,
$\alpha\in B_n$, where $\theta$ is the induced permutation of $\alpha$.
The semidirect product $B_n\ltimes G^n$ is a Garside group.
Note that if $G$ is the $m$-braid group $B_m$,
$B_n\ltimes (B_m)^n$ consists of reducible braids.
\end{exmp}

The second example is similar to the wreath product, which is
a semidirect product of a cartesian product of a group and a subgroup of
a symmetric group, where the action is the coordinate change.
Since Garside groups are torsion-free, we use
a braid group instead of a subgroup of a symmetric group.

\begin{lem}\label{thm:infsup}
Let $G\ltimes H$ be the semidirect product of Garside groups $G$ and $H$
via a homomorphism $\rho\colon G\to\Aut(H^+)$.
Let $\Delta_G$ and $\Delta_H$ be Garside elements of\/ $G$ and $H$ such that
$\Delta_H$ is fixed under the action of\/ $G$.
For $(g,h)\in G\ltimes H$, the following hold, where
$(\Delta_G,\Delta_H)$ is the Garside element of $G\ltimes H$ used
for {\upshape(3), (4)} and {\upshape(5)}.
\begin{enumerate}
\item $(g,h)^{-1}=(g^{-1}, (h^{-1})^{g^{-1}})$.
\item $(\Delta_G,\Delta_H)^k=(\Delta_G^k,\Delta_H^k)$ for all $k\in\Z$.
\item $\tau((g,h))=(\tau(g),\tau(h^{\Delta_G}))$.
\item $\inf((g,h))=\min\{\inf(g),\inf(h)\}$.
\item $\sup((g,h))=\max\{\sup(g),\sup(h)\}$.
\end{enumerate}
\end{lem}

\begin{pf}
(1) and (2) are obvious.

\smallskip
(3)\ \
$\tau((g,h))
=(\Delta_G^{-1},\Delta_H^{-1})(g,h)(\Delta_G,\Delta_H)
=(\Delta_G^{-1} g\Delta_G,\Delta_H^{-1} h^{\Delta_G} \Delta_H)
=(\tau(g),\tau(h^{\Delta_G}))$.

\smallskip
(4)\ \
By Proposition~\ref{thm:SemidirectProperties}(4),
$(\Delta_G^k,\Delta_H^k)\le_R(g,h)$ if and only if $\Delta_G^k\le_R g$
and $\Delta_H^k\le_R h$.
Since $(\Delta_G,\Delta_H)^k=(\Delta_G^k,\Delta_H^k)$,
$\inf(g,h)=\min\{\inf(g),\inf(h)\}$.

\smallskip
(5) can be proved similarly to (4).
\qed\end{pf}

We close this section showing that
if\/ $\Delta_G$ and $\Delta_H$ are the minimal Garside
elements of $G$ and $H$, respectively,
then $(\Delta_G,\Delta_H)$ is the
minimal Garside element of\/ $G\ltimes H$.
Because we don't have a proof using only the techniques
developed in this paper, we use the
the characterization of Garside elements
by Dehornoy~\cite[Proposition 1.10]{Deh02}.

Let us say a subset $S$ of a Garside monoid
\emph{LC-closed} if the following holds.
\begin{itemize}
\item If $a,b\in S$, then
$a\wedge_L b$, $a\wedge_R b$, $a\vee_L b$,
$a\vee_R b$, $a^{-1}(a\vee_L b)$, $(a\vee_R b)a^{-1}\in S$.
\item If $a\in S$ and $b\le_L a$ or $b\le_R a$, then $b\in S$.
\end{itemize}
(LC stands for lattice operations and complements.
The elements $a^{-1}(a\vee_L b)$ and $(a\vee_R b)a^{-1}$ are called
the right complement of $a$ in $b$ and the left complement of $a$ in $b$
respectively~\cite{Deh02}.)
For a subset $A$ of a Garside monoid, let \emph{LC-closure}
of $A$ be the smallest LC-closed set containing $A$.
Then the following holds.

\begin{itemize}
\item $\Delta$ is a Garside element and
$\D$ is the set of simple elements corresponding to $\Delta$ if and only if
$\D$ is an LC-closed finite generating set of
the Garside monoid and $\Delta$ is the right lcm of the elements
of $\D$~\cite[Proposition 1.10]{Deh02}.

\item $\Delta$ is the minimal Garside element and
$\D$ is the set of simple elements  corresponding to $\Delta$ if and only if
$\D$ is the LC-closure of the set of atoms
and $\Delta$ is the right lcm of all elements
of $\D$~\cite[discussion after Definition 1.11 in page 275]{Deh02}.
\end{itemize}

\begin{lem}
If\/ $\Delta_G$ and $\Delta_H$ are the minimal Garside
elements of\/ $G$ and $H$, respectively,
then $(\Delta_G,\Delta_H)$ is the
minimal Garside element of\/ $G\ltimes H$.
\end{lem}

\begin{pf}
Let $\Delta_G$, $\Delta_H$ and $\Delta$
be the minimal Garside elements of
$G$, $H$ and $G\ltimes H$, respectively.
We show that $(\Delta_G,\Delta_H)$ is equal to $\Delta$.

Because $(\Delta_G,\Delta_H)$ is a Garside element of $G\ltimes H$
and $\Delta$ is the minimal Garside element of $G\ltimes H$,
we have $\Delta\le_R(\Delta_G,\Delta_H)$.

Let $A=\{(a,1):\mbox{$a$ is an atom of $G^+$}\}$
and $B=\{(1,b):\mbox{$b$ is an atom of $H^+$}\}$.
It is easy to see that the LC-closure of $A$ is
$\{(s,1): s\le\Delta_G\}$, observing that for $a_1,a_2\in G^+$
$(a_1,1)\wedge_L(a_2,1)=(a_1\wedge_L a_2,1)$ and
the same is true for $\vee_L$, $\wedge_R$ and $\vee_R$.
In particular, the LC-closure of $A$ contains $(\Delta_G,1)$.
Similarly, the LC-closure of $B$ contains $(1,\Delta_H)$.
Therefore the LC-closure of $A\cup B$ contains
$(\Delta_G,1)\vee_R (1,\Delta_H)=(\Delta_G,\Delta_H)$.
Because $A\cup B$ is the set of atoms of $G\ltimes H$,
$\Delta$ is the right lcm of the element of
the LC-closure of $A\cup B$.
Hence $(\Delta_G,\Delta_H)\le_R\Delta$.

Consequently, $(\Delta_G,\Delta_H)$ is equal to $\Delta$.
\qed\end{pf}

\section{The product $G(n)=\Z\ltimes G^n$}

\begin{defn}
For a Garside group $G$, we define $G(n)$
as the semidirect product $\Z\ltimes G^n$,
where $\Z=\langle \delta\rangle$ acts on the cartesian product $G^n$
by $(g_1,\ldots,g_n)^\delta=(g_n,g_1,\ldots,g_{n-1})$.
The element $(\delta^k,(g_1,\ldots,g_n))\in G(n)$
is denoted $\delta^k(g_1,\ldots,g_n)$.
\end{defn}

The infinite cyclic group $\Z=\langle \delta\rangle$ is a Garside group
with minimal Garside element $\delta$.
(Observe that $\Z$ is isomorphic to the 2-strand braid group $B_2$.)
Therefore, $G(n)$ is a Garside group and
if $\Delta$ is a Garside element of $G$, then
$\delta(\Delta,\ldots,\Delta)$ is a Garside element of $G(n)$.
Throughout the entire section, the Garside elements
$\delta$, $\Delta$ and $\delta(\Delta,\ldots,\Delta)$
are used for the respective groups.

\begin{lem}\label{thm:inf-Gn}
In $G(n)$, the following hold.
\begin{enumerate}
\item $\inf(\delta^k(g_1,\ldots,g_n))=\min\{k,\inf(g_1),\ldots,\inf(g_n)\}$.
\item $\sup(\delta^k(g_1,\ldots,g_n))=\max\{k,\sup(g_1),\ldots,\sup(g_n)\}$.
\item $\tau(\delta^k(g,\ldots,g))=\delta^k(\tau(g),\ldots,\tau(g))$.
\item
$\c(\delta^k(g,\ldots,g))=\biggl\{\begin{array}{ll}
\delta^k(\c(g),\ldots,\c(g)) &\mbox{if\/ $k\ge\inf(g)$},\\
\delta^k(\tau(g),\ldots,\tau(g)) &\mbox{if\/ $k<\inf(g)$.}
\end{array}$
\item
$\d(\delta^k(g,\ldots,g))=\biggl\{\begin{array}{ll}
\delta^k(\d(g),\ldots,\d(g)) &\mbox{if\/ $k\le\sup(g)$},\\
\delta^k(g,\ldots,g) &\mbox{if\/ $k>\sup(g)$.}
\end{array}$
\end{enumerate}
\end{lem}

\begin{pf}
(1), (2) and (3) are immediate consequences of Lemma~\ref{thm:infsup}.
We prove (4) for the case $k=\inf(g)$.
The other cases of (4) and (5) can be proved similarly.
Let $\alpha=\delta^k(g,\ldots,g)$ and $g=\Delta^k sa$, where $s=\Delta\wedgeL (sa)$.
Let $\DD=\delta(\Delta,\ldots,\Delta)$.
Then, $\inf(\alpha)=k$ by (1) and
$$
\alpha=\delta^k(g,\ldots,g)
=\delta^k(\Delta^k sa,\ldots,\Delta^k sa)
=\DD^k(sa,\ldots,sa).
$$
Since $\DD\wedgeL(sa,\ldots,sa)=(s,\ldots,s)$
by Proposition~\ref{thm:SemidirectProperties}~(3),
\begin{eqnarray*}
\c(\alpha)
&=&\DD^k(a,\ldots,a)\tau^{-k}((s,\ldots,s))
=\DD^k(a,\ldots,a)(\tau^{-k}(s),\ldots,\tau^{-k}(s))\\
&=&\DD^k(a\tau^{-k}(s),\ldots,a\tau^{-k}(s))
= \delta^k(\Delta^ka\tau^{-k}(s),\ldots,\Delta^ka\tau^{-k}(s))\\
&=&\delta^k(\c(g),\ldots,\c(g)).\qed
\end{eqnarray*}
\end{pf}

\begin{cor}\label{thm:sss-gn}
Let $\alpha=\delta^k(g,\ldots,g)\in G(n)$.
\begin{enumerate}
\item The ultra summit set of\/ $\alpha$ contains an element
of the form $\delta^m(h,\ldots,h)$.
\item If $g$ is contained in its super summit set, then so is $\alpha$.
\item If $g$ is contained in its ultra summit set, then so is $\alpha$.
\end{enumerate}
\end{cor}

\begin{pf}
(1)\ \
There exists $l_1,l_2\ge 0$ such that $\c^{l_1}\d^{l_2}(\alpha)\in[\alpha]^U$.
By Lemma~\ref{thm:inf-Gn} (4) and (5), $\c^{l_1}\d^{l_2}(\alpha)$
is of the form $\delta^m(h,\ldots,h)$.

\smallskip
(2)\ \
Note that for $h\in[g]^S$ and $l\ge 1$,
$\inf(h)=\inf(\c^l(h))=\inf(\d^l(h))=\inf(\tau^l(h))$.

If $k\ge\inf(g)$, then $k\ge\inf(\c^l(g))$ for all $l\ge 1$, hence
$\c^l(\alpha)=\delta^k(\c^l(g),\ldots,\c^l(g))$
by Lemma~\ref{thm:inf-Gn}~(4).
Since $\inf(\c^l(\alpha))=\min\{k,\inf(\c^l(g))\}
=\min\{k,\inf(g)\}=\inf(\alpha)$ for all $l\ge 1$,
we obtain $\inf(\alpha)=\infs(\alpha)$.

Similarly, if $k<\inf(g)$, then $k<\inf(\c^l(g))$ for all $l\ge 1$, hence
$\c^l(\alpha)=\delta^k(\tau^l(g),\ldots,\tau^l(g))$
by Lemma~\ref{thm:inf-Gn}~(4).
Since $\inf(\c^l(\alpha))=\min\{k,\inf(\tau^l(g))\}
=\min\{k,\inf(g)\}=\inf(\alpha)$ for all $l\ge 1$,
we obtain $\inf(\alpha)=\infs(\alpha)$.

In both cases, we have $\inf(\alpha)=\infs(\alpha)$.
Similarly, we obtain that $\sup(\alpha)=\sups(\alpha)$,
hence $\alpha$ is contained in its super summit set.

\smallskip
(3)\ \
If $k<\inf(g)$, then for any $l\ge 1$,
$\c^l(\alpha)=\delta^k(\tau^l(g),\ldots,\tau^l(g))$ as in (2).
Note that $\Delta^e$ is central for some $e\ge 1$.
Then $\tau^e(g)=g$, hence $\c^e(\alpha)=\alpha$ and $\alpha\in[\alpha]^U$.

If $k\ge\inf(g)$, then for any $l\ge 1$,
$\c^l(\alpha)=\delta^k(\c^l(g),\ldots,\c^l(g))$ as in (2).
Because $g$ is contained in its ultra summit set,
$\c^m(g)=g$ for some $m\ge 1$,
hence $\c^m(\alpha)=\alpha$ and $\alpha\in[\alpha]^U$.
\qed\end{pf}

\begin{lem}\label{thm:conj-gn}
If\/ $k\equiv 1\pmod n$, then
$\delta^k(g_1,\ldots,g_n)$ and $\delta^k(h_1,\ldots,h_n)$ are conjugate in $G(n)$
if and only if\/ $g_1\cdots g_n$ and $h_1\cdots h_n$ are conjugate in $G$.
\end{lem}

\begin{pf}
Note that $\delta^k(g_1,\ldots,g_m,g_{m+1},1,\ldots,1)$
and $\delta^k(g_1,\ldots,g_{m-1},g_mg_{m+1},1,1,\ldots,1)$
are conjugate for $1\le m\le n-1$.
(The latter can be obtained from the former
by conjugating on the right by the element
$(1,\ldots,1,g_{m+1},1,\ldots,1)$,
where the nontrivial coordinate is in the $m$th position.)
Therefore,
$$
\delta^k(g_1,\ldots,g_n)
\quad\mbox{is conjugate to}\quad \delta^k(g_1\cdots g_n,1,\ldots,1).
\eqno(\ast)
$$

Assume that $g=g_1\cdots g_n$ is conjugate to $h=h_1\cdots h_n$ in $G$.
Let $g=x^{-1}hx$. Since
$$
(x^{-1},\ldots,x^{-1})\delta^k(h,1,\ldots,1)(x,\ldots,x)=\delta^k(g,1,\ldots,1),
$$
$\delta^k(g,1,\ldots,1)$ is conjugate to $\delta^k(h,1,\ldots,1)$, and thus
$\delta^k(g_1,\ldots,g_n)$ and
$\delta^k(h_1,\ldots,h_n)$ are conjugate by~($\ast$).

Conversely, assume that $\delta^k(g_1,\ldots,g_n)$ is conjugate to
$\delta^k(h_1,\ldots,h_n)$ in $G(n)$.
Let $g=g_1\cdots g_n$ and $h=h_1\cdots h_n$.
Then $\delta^k(g,1,\ldots,1)$ is conjugate to
$\delta^k(h,1,\ldots,1)$ by~($\ast$).
Let $\delta^k(g,1,\ldots,1)=\gamma^{-1}\delta^k(h,1,\ldots,1)\gamma$ for
some $\gamma=\delta^m(x_1,\ldots,x_n)$. Assume $m\not\equiv 0\pmod n$.
Since
\begin{eqnarray*}
\lefteqn{\delta^k(g,1,\ldots,1)=\gamma^{-1}\delta^k(h,1,\ldots,1)\gamma}\\
&=&(x_1^{-1},\ldots,x_n^{-1})
    \delta^{-m}\delta^k(h,1,\ldots,1)
    \delta^m(x_1,\ldots,x_n)\\
&=&\delta^k(x_n^{-1},x_1^{-1},\ldots,x_{n-1}^{-1})
  (\underbrace{1,\ldots,1}_m,h,
  \underbrace{1,\ldots,1}_{n-m-1})
  (x_1,\ldots,x_n)\\
&=&\delta^k(x_n^{-1} x_1, x_1^{-1} x_2,\ldots,x_{m-1}^{-1} x_m,
   x_m^{-1} h x_{m+1}, x_{m+1}^{-1}x_{m+2},\ldots, x_{n-1}^{-1}x_n),
\end{eqnarray*}
$x_i=x_{i+1}$ for $i\ne m,n$; $g=x_n^{-1} x_1$;  $x_m=hx_{m+1}$.
(Here, $x_k$ denotes $x_j$
such that $1\le j\le n$ and $k\equiv j\pmod n$.)
Since
$x_ng=x_1=x_2=\cdots=x_m=hx_{m+1}=hx_{m+2}=\cdots =hx_n$,
$g=x_n^{-1} h x_n$. Therefore $g$ is conjugate to $h$.
The case $m\equiv 0\pmod n$ can be proved similarly.
\qed\end{pf}

\section{$\infs$ and $\sups$ of powers}

\begin{thm}\label{thm:infPower}
Let $G$ be a Garside group. For $g\in G$ and $n\ge 1$,
\begin{enumerate}
\item $\infs(g)  \le \frac{\infs(g^n)}n < \infs(g)+1$,
\item $\sups(g)-1 <  \frac{\sups(g^n)}n \le \sups(g)$,
\item $\lens(g)-2 <  \frac{\lens(g^n)}n \le \lens(g)$.
\end{enumerate}
\end{thm}

\begin{pf}
We prove only (1), because (2) follows from (1) and Lemma~\ref{thm:dual},
and (3) follows from (1) and (2).
It is obvious that $\infs(g)\le\frac{\infs(g^n)}n$.
Since $\infs(g)$ depends only on the conjugacy class of $g$,
we may assume that $g$ is contained in its super summit set.
Let $\infs(g)=\inf(g)=r$.
Choose an integer $k$ such that $k>r$ and $k\equiv 1\pmod n$.
Let $\alpha=\delta^k(g,\ldots,g)\in G(n)$.
Then $\alpha$ is contained in its super summit set
by Corollary~\ref{thm:sss-gn}~(2).
By Lemma~\ref{thm:inf-Gn}~(1),
$$
\infs(\alpha)=\inf(\alpha)=\min\{k,\inf(g)\}=\min\{k,r\}=r.\eqno(\ast)
$$
Assume that $\inf_s(g^n)\ge n(\infs(g)+1)=n(r+1)$.
Then $g^n$ is conjugate to $\Delta^{n(r+1)}a$ for some $a\in G^+$.
Since $k\equiv 1\pmod n$,
$\alpha=\delta^k(g,\ldots,g)$ and
$\beta=\delta^k(\Delta^{r+1},\ldots,\Delta^{r+1},\Delta^{r+1}a)$
are conjugate by Lemma~\ref{thm:conj-gn}.
By Lemma~\ref{thm:inf-Gn}~(1),
$$\infs(\alpha)\ge\inf(\beta)=\min\{k,r+1\}=r+1,$$
which contradicts $(\ast)$.
Therefore $\infs(g^n)< n(\infs(g)+1)$.
\qed\end{pf}

We remark that the inequalities in Theorem~\ref{thm:infPower}
do not hold for $\inf$, $\sup$ or $\len$.
For example, consider the 3-braid group
$B_3=\langle\sigma_1,\sigma_2\mid
\sigma_1\sigma_2\sigma_1=\sigma_2\sigma_1\sigma_2\rangle$.
It is a Garside group with minimal Garside element
$\Delta=\sigma_1\sigma_2\sigma_1$.
Let $g=\sigma_1^{-(k+1)}\Delta\sigma_1^{k+1}=\sigma_1^{-k}\sigma_2\sigma_1^{k+2}$.
Then $\inf(g)=-k$ and $\infs(g)=1$, and for any $m\ge 1$,
$\inf(g^{2m})=2m$
because $g^{2m}=\Delta^{2m}$ and $\Delta^2$ is central.
Therefore,
$\frac{\inf(g^{2m})}{2m}=1\ge -k+1=\inf(g)+1$ for any $k\ge 0$ and $m\ge 1$.

Let $\lfloor x\rfloor$ denote the largest integer less than or
equal to $x$ and $\lceil x\rceil$ the smallest integer greater
than or equal to $x$.

\begin{cor}\label{thm:infUnique}
Let $G$ be a Garside group, $g\in G$ and $n\ge 1$.
Then $\inf_s(g)$ and\/ $\sup_s(g)$ are uniquely determined
by $(n, \inf_s(g^n), \sup_s(g^n))$ as follows.
$$
\inf{\!}_s(g)=\left\lfloor\frac{\infs(g^n)}{n}\right\rfloor\quad\mbox{and}\quad
\sup{\!}_s(g)=\left\lceil\frac{\sups(g^n)}{n}\right\rceil.
$$
In particular, for any $n\ge 1$, $\inf_s(g)\ge 0$ if and only if\/
$\infs(g^n)\ge 0$ and $\sup_s(g)\le 0$ if and only if\/
$\sup_s(g^n)\le 0$.
\end{cor}

\begin{cor}
If\/ $g^n$ is conjugate to $g^{-n}$ for some $n\ne 0$, then
$\inf_s(g)=-\sup_s(g)$.
\end{cor}

\begin{pf}
We may assume that $n\ge 1$. Since $g^n$ is conjugate to
$g^{-n}$, $\infs(g^n)=\infs(g^{-n})=-\sups(g^n)$. By
Corollary~\ref{thm:infUnique}, we are done.
\qed\end{pf}

\section{Garside groups are strongly translation discrete}
\label{sec:trans}

Generalizing Bestvina's work on the Artin groups of finite type,
Charney, Meier and Whittlesey~\cite{CMW04}
proved that Garside groups are translation discrete,
assuming the existence of a tame Garside element.
A Garside element $\Delta$ is said to be \emph{tame}
if there exists a constant $c>0$ such that
$\Vert\Delta^n\Vert \le cn$ for all $n\ge 1$.
This property has recently been explored
in~\cite{Sib04}, where general criteria are established
that imply tameness and where it is conjectured that all Garside
monoids have a tame Garside element.
In this section, we prove the main theorem of this paper
that all Garside groups are strongly translation
discrete without any assumption.

Let $G$ be a Garside group and $\D$ the set of simple elements.
Abusing the notation, let $|g|_\D$ denote
the shortest word length of $g\in G$ in terms of $\D\cup\D^{-1}$,
in other words, $|g|_\D$ is $|g|_{\D\cup\D^{-1}}$ in the
notation of~\S1.  Similarly, let $t_\D(g)$ denote
the translation number of $g\in G$ with respect to $\D\cup \D^{-1}$.

The shortest word length $|g|_\D$ is also known as \emph{geodesic length}.
It is known that~\cite{Cha95}
$$
|g|_\D=\left\{\begin{array}{l}
\sup(g) \quad\mbox{if $\inf(g)\ge 0$};\\
-\inf(g)\quad\mbox{if $\sup(g)\le 0$};\\
\len(g) \quad\mbox{if $\inf(g)<0<\sup(g)$.}
\end{array}\right.
$$

\begin{thm}\label{thm:TranUltra}
Let $G$ be a Garside group and $\D$ the set of simple elements.
If\/ $g\in G$ is contained in its super summit set,
then for any $n\ge 1$
$$
|g|_\D-2< \frac{|g^n|_\D}n\le |g|_\D.
$$
In particular, $|g|_\D-2\le t_\D(g)\le |g|_\D$.
\end{thm}

\begin{pf}
It is obvious that $\frac{|g^n|_\D}n\le |g|_\D$.
Since $g$ is  in its super summit set,
$\inf(g)=\infs(g)$, $\sup(g)=\sups(g)$
and $\len(g)=\lens(g)$.
If $\inf(g)\ge 0$, then $\inf(g^n)\ge n\inf(g)\ge 0$, and thus
$$
\frac{|g^n|_\D}n
=\frac{\sup(g^n)}n
\ge \frac{\sups(g^n)}n
> \sups(g)-1
=\sup(g)-1
=|g|_\D-1.
$$
If $\sup(g)\le 0$, then $\sup(g^n)\le n\sup(g)\le 0$, and thus
$$
\frac{|g^n|_\D}n
=\frac{-\inf(g^n)}n
\ge\frac{-\infs(g^n)}n
> -\infs(g)-1
=-\inf(g)-1
=|g|_\D-1.
$$
If $\inf(g)<0<\sup(g)$, then $\infs(g^n)<0<\sups(g^n)$
by Corollary~\ref{thm:infUnique}.
Therefore, $\inf(g^n)<0<\sup(g^n)$ and
$$\frac{|g^n|_\D}n
=\frac{\len(g^n)}n\ge \frac{\lens(g^n)}n
> \lens(g)-2
= \len(g)-2
=|g|_\D-2.
$$
In all cases, $\frac{|g^n|_\D}n>|g|_\D-2$.
Therefore, $|g|_\D-2< \frac{|g^n|_\D}n\le |g|_\D$ as desired.
By taking $n\to\infty$, $|g|_\D-2\le t_\D(g)\le |g|_\D$.
\qed\end{pf}

\begin{thm}\label{thm:main}
All Garside groups are strongly translation discrete.
\end{thm}

\begin{pf}
Let $G$ be a Garside group and $\D$ the set of simple elements.
Since Garside groups are biautomatic by Dehornoy~\cite{Deh02}
and biautomatic groups are translation separable
by Gersten and Short~\cite{GS91},
$G$ is translation separable.
Therefore, it suffices to show that for any real number $r$,
there are only finitely many conjugacy classes $[g]$ such that
$t_\D(g)\le r$.
By Theorem~\ref{thm:TranUltra},
$|g|_\D\le t_\D(g)+2$ if $g\in[g]^S$.
Since $t_\D(g)$ depends only on the conjugacy class of $g$,
$$
\{ [g]: t_\D(g)\le r\}
=\{ [g]: t_\D(g)\le r,\ g\in[g]^S\}
\subset \{ [g]: |g|_\D \le r+2\}.
$$
Since there are only finitely many elements $g$ such that $|g|_\D \le r+2$,
the set $\{[g]: t_\D(g)\le r\}$ is finite.
\qed\end{pf}

\begin{cor}
Let $G$ be a Garside group.
\begin{enumerate}
\item Every solvable subgroup of $G$ is
finitely generated and virtually abelian.
\item $G$ cannot contain subgroups isomorphic to
the additive group of rational numbers
or the group of $p$-adic fractions
$\Q_p=\{k/p^l\mid k\in\Z,\ l\in\N\}$.
\item For any finite set $X$ of semigroup generators for $G$,
$t_X(G)$ is a closed discrete set.
\end{enumerate}
\end{cor}

\begin{pf}
The proof of (1) is the same as for Bestvina's Corollaries 4.2 and 4.4
in~\cite{Be99}.
(2) is a property of translation discrete groups,
and (3) is a property of strongly translation discrete groups.
See~\cite{Kap97} for (2) and (3).
\qed\end{pf}

\section{The root problem in a Garside group}
\label{sec:root}

Generalizing the work of Sty\v shnev~\cite{Sty78} on the braid groups,
Sibert~\cite{Sib02} showed that the root problem
is solvable for Garside groups
under a certain condition, called the finiteness of
positive conjugacy classes.
Two elements $a$ and $b$ of a Garside monoid $G^+$ are said to be
\emph{positively conjugate} if there exists $c\in G^+$ satisfying $ac=cb$.
Note that two elements of $G^+$ are postively conjugate
if and only if they are conjugate in $G$,
because for any element $g$ of $G$, there exists an integer $m$ such that
$\Delta^m$ is central and $\Delta^m g\in G^+$.
The assumption of Sibert is that
the positive conjugacy classes are all finite.
It is known that if a Garside monoid has a tame Garside element,
then it satisfies the finiteness condition of
positive conjugacy class~\cite{Sib04}.

In this section, we present a new solution to the root problem
in Garside groups and show that the root problem is solvable for
Garside groups without any assumption.

\begin{thm}
Let $G$ be a Garside group and $g\in G$.
\begin{enumerate}
\item There are only finitely many $n\in\N$ such that $g$ has an $n$th root.
\item For each $n\in\N$, there are only finitely
many conjugacy classes of $n$th roots of\/ $g$.
\end{enumerate}
\end{thm}

\begin{pf}
(1) is a direct consequence of the obvious property $t_X(g^n)=n\,t_X(g)$
of translation numbers, together with the definition
of translation discreteness.
See~\cite[page 1853]{Kap97}.
(2) is a direct consequence of (1), together with the definition
of strong translation discreteness.
\qed\end{pf}

Recall that the $n$th root of a braid is unique up to conjugacy
by Gonz\'alez-Meneses~\cite{Gon04}.
However, this is not true for Garside groups.
For example, the group
$$
G=\langle
x_1,\ldots,x_m:x_1^{p_1}=x_2^{p_2}=\cdots=x_m^{p_m}\rangle,\qquad
p_i\ge 2,
$$
is a Garside group with a Garside element
$\Delta=x_1^{p_1}$. See~\cite[Example 4]{DP99}.
If the exponents $p_1,\ldots,p_m$ have a
common divisor $n\ge 2$, then there are at least $m$ conjugacy
classes of $n$th roots of $\Delta$.

\begin{thm}\label{thm:rootprob_Gn}
Let\/ $G$ be a Garside group and $g\in G$.
The element $g$ has an $n$th root if and
only if the ultra summit set of\/ $\delta(g,1,\ldots,1)\in G(n)$
contains an element of the form $\delta(h,\ldots,h)$. In this case,
if\/ $\delta(g,1,\ldots,1)=\gamma^{-1}\delta(h,\ldots,h)\gamma$ for
some $\gamma=\delta^k(x_1,\ldots,x_n)$,
then $g=(x_n^{-1} hx_n)^n$.
\end{thm}

\begin{pf}
The first statement is immediate from Corollary~\ref{thm:sss-gn}
and Lemma~\ref{thm:conj-gn}.
Suppose
$\delta(g,1,\ldots,1)=\gamma^{-1}\delta(h,\ldots,h)\gamma$ for
$\gamma=\delta^k(x_1,\ldots,x_n)$. Since
$$
\delta(g,1,\ldots,1)=\gamma^{-1}\delta(h,\ldots,h)\gamma =
\delta(x_n^{-1} h x_1, x_1^{-1} hx_2,\ldots x_{n-1}^{-1}hx_n),
$$
$g=x_n^{-1} hx_1$ and $x_i=hx_{i+1}$ for $i=1,\ldots,n-1$.
Therefore, $x_1=hx_2=h^2x_3=\cdots=h^{n-1}x_n$,
and thus
$g=x_n^{-1}hx_1=x_n^{-1}h(h^{n-1}x_n)=x_n^{-1} h^n x_n
=(x_n^{-1}hx_n)^n$.
\qed\end{pf}

\begin{cor}\label{thm:rootprob}
The root problem is solvable for the Garside groups.
\end{cor}

We remark that Theorem~\ref{thm:rootprob_Gn}
does not require the condition of
finite positive conjugacy class of~\cite{Sib02}.
The following algorithm solves the root problem.
It seems more efficient than
the algorithm of Sty\v shnev~\cite{Sty78}
and Sibert~\cite{Sib02},
because it uses the conjugacy algorithm in
$G(n)$ instead of an exhaustive search in $G^+$.

\smallskip
\noindent{\bf Root Extraction Algorithm}.\ \ Given a Garside group
$G$, $g\in G$ and an integer $n\ge 2$:
\begin{itemize}
\item[\bf 1.] Set $\alpha=\delta(g,1,\ldots,1)\in G(n)$ and
compute the ultra summit set $[\alpha]^U$ of $\alpha$. \item[\bf
2.] For each $\beta\in[\alpha]^U$, test whether $\beta$ is of the
form $\delta(h,\ldots,h)$ for some $h\in G$. \item[]\mbox{}\qquad
If so, compute $\gamma=\delta^k(x_1,\ldots,x_n)$ such that
$\alpha=\gamma^{-1}\beta\gamma$. Then $g=(x_n^{-1}hx_n)^n$.
\item[\bf 4.] If there is no element in $[\alpha]^U$ of the form
$\delta(h,\ldots,h)$, then $g$ has no $n$th root.
\end{itemize}

\section*{Acknowledgement}
The author thanks the referees
for their careful reading and useful suggestions.
He also thanks the Korea Institute for Advanced Study
(KIAS). Part of the work was done during the stay of the
author at this institute.

\end{document}